\documentclass[11pt, english]{article}
\usepackage[
    a4paper
    ]{geometry}
\usepackage[utf8]{inputenc}
\usepackage{lipsum}
\usepackage{babel, amsmath, amsthm, amssymb, amsfonts, mathtools, centernot, enumitem, mathrsfs}
\usepackage{thmtools, thm-restate}
\usepackage{tikz-cd}
\usepackage{tikz-3dplot}
\usepackage{caption}
\usepackage{intcalc}
\usepackage{stmaryrd}
\usepackage{multicol}
\usepackage{cite}
\usepackage{hyperref}
\usepackage{caption}
\usepackage{subcaption}
\usepackage[capitalise]{cleveref}
\usepackage[toc,page]{appendix}

\usepackage[toc,page]{appendix}
\usepackage{fancyhdr}
\usepackage{todonotes}
\newcommand\tab[1][1cm]{\hspace*{1}}

\newtheorem{theorem}{Theorem}[section]
\declaretheorem[name=Theorem,sibling=theorem]{restate-thm}

\newtheorem{corollary}{Corollary}[theorem]
\newtheorem{proposition}[theorem]{Proposition}
\newtheorem{lemma}[theorem]{Lemma}
\theoremstyle{definition}
\newtheorem{definition}[theorem]{Definition}
\newtheorem{example}[theorem]{Example}
\newtheorem*{problem}{Question}
\newtheorem{mainthm}{Theorem}

\newtheoremstyle{TheoremNum}
        {\topsep}{\topsep}              
        {\itshape}                      
        {}                              
        {\bfseries}                     
        {.}                             
        { }                             
        {\thmname{#1}\thmnote{ \bfseries #3}}
    \theoremstyle{TheoremNum}

\renewcommand{\P}{\mathcal P}
\newcommand{\E}{\mathcal E}
\newcommand{\M}{\mathcal M}
\newcommand{\eff}{{}^\perp \mathcal P}
\DeclareMathOperator{\Hom}{Hom}
\DeclareMathOperator{\cogen}{cogen}
\DeclareMathOperator{\Ext}{Ext}

\DeclareMathOperator{\End}{End}

\DeclareMathOperator{\Image}{Im}

\DeclareMathOperator{\domdim}{domdim}
\DeclareMathOperator{\gldim}{gl.dim}

\DeclareMathOperator{\id}{id}
\DeclareMathOperator{\add}{add}
\DeclareMathOperator{\Mod}{Mod}
\DeclareMathOperator{\adm}{adm}
\def\mod{\operatorname{mod}}

\newcommand{\pac}{Morita--Tachikawa category}

\newsavebox{\pullback}
\sbox\pullback{%
	\begin{tikzpicture}%
	\draw (0,0) -- (1ex,0ex);%
	\draw (1ex,0ex) -- (1ex,1ex);%
	\end{tikzpicture}
}

\begin{document}
\title{Higher Auslander--Solberg correspondence for exact categories}
\author{Jacob Fjeld Grevstad}
\date{2022}

\maketitle
\begin{abstract}
  \noindent We introduce the concept of an $n$-minimal Auslander--Gorenstein category and $n$-precluster tilting subcategory. With this, we create an analog of the higher Auslander--Solberg correspondence (Iyama--Solberg 2018) for exact categories. Our approach is based on the recent generalization of the (higher) Auslander correspondence to exact categories (Henrard--Kvamme--Roosmalen 2020,  Ebrahimi--Nasr-Isfahani 2021).
\end{abstract}

\newpage
\tableofcontents
\newpage

\section*{Introduction}
\addcontentsline{toc}{section}{\numberline{}Introduction}

\addcontentsline{toc}{subsection}{\numberline{}History}

In the early '70s Auslander showed a correspondence between representation finite algebras, and algebras $\Gamma$ satisfying 
$$\gldim \Gamma \leq 2 \leq \domdim\Gamma$$ 
up to Morita equivalence \cite{Aus71}. Historically, this was the starting point of Auslander--Reiten theory.

In more recent years the Auslander correspondence has been generalized in several ways. In 2007, Iyama established a correspondence between $n$-cluster tilting modules and algebras satisfying 
$$\gldim \Gamma \leq n+1 \leq \domdim\Gamma,$$ 
known as the higher Auslander correspondence\cite{Iyama07}. In the same way the Auslander-correspondence is significant for Auslander--Reiten theory, $n$-cluster tilting modules form an important part of higher Auslander--Reiten theory. One way this manifests is that $n$-cluster tilting modules satisfy $\tau_n(M) \in \add M$ and $\tau^-_n(M)\in \add M$, where $\tau_n$ and $\tau^-_n$ are the higher Auslander--Retien translates.

Inspired by this property, Iyama and Solberg defined $n$-precluster tilting modules\cite{IS18} to be $n$-rigid generator-cogenerators which satisfied $\tau_n(M) \in \add M$ and $\tau^-_n(M)\in \add M$. They also established a higher Auslander--Solberg correspondence between $n$-precluster tilting modules and Gorenstein algebras satisfying 
$$\id_\Gamma \Gamma \leq n+1 \leq \domdim\Gamma,$$ 
replacing finite global dimension with being Gorenstein. The case of $n=1$ was proved by Auslander and Solberg in 1993\cite{AS93gorenstein}, hence the name.

In 2020 Henrard, Kvamme and Roosmalen introduced an analog of the Auslander correspondence for exact categories. Instead of taking the endomorphism ring of an additive generator for a representation finite algebra, they take an arbitrary exact category an consider the category of admissibly presented modules. This creates a correspondence between exact categories and what they call Auslander exact categories.

The following year Ebrahimi and Nasr-Isfahani showed that the construction of Henrard--Kvamme--Roosmalen also gives a correspondence between $n$-cluster tilting subcategories and $n$-Auslander exact categories.

In this paper, in much the same way, we introduce the concept of an $n$-precluster tilting subcategory of an exact category and the concept of an $n$-minimal Auslander--Gorenstein exact category. We also show that the construction by Henrard--Kvamme--Roosmalen gives a bijection between these.

\subsubsection*{Structure and main theorems}
\addcontentsline{toc}{subsection}{\numberline{}Structure and main theorems}

The goal of this paper is to generalize the correspondence in \cite{henrard2020auslanders} to create an analog of the higher Auslander--Solberg correspondence for exact categories. In \cref{sec:preliminaries} we go over some basic results and definition relating to exact categories, and in \cref{sec:admissible_modules} we recall some properties of the category of admissible modules from \cite{henrard2020auslanders}.

The main theory of this paper is in \cref{sec:Auslander-Solber-correspondence} where we discuss a few different related correspondences.

The Morita--Tachikawa correspondence \cite{GK15, Ringel07} says that there is a bijection between generator-cogenerators and algebras with dominant dimension at least 2. Since we will be working in exact categories which don't necessarily have enough injectives, the concept of dominant dimension is not directly useful. Therefore we define Morita--Tachikawa categories (\cref{def:pac}) as an exact analog to algebras with dominant dimension at least 2. We do this in an analogues way to how Henrard--Kvamme--Roosmalen define Auslander exact categories \cite[Def~1.3]{henrard2020auslanders}. This gives us our first main theorem.

\begin{mainthm}[\cref{thm:bijection_partial_auslander}]
  There is a bijection between equivalence classes of contravariantly finite, generating-cogenerating subcategories and equivalence classes of Morita--Tachikawa categories.
\end{mainthm}

The proof of this theorem is implicit in \cite{ebrahimi2021higher}, but they do not state it explicitly.

Another generalization of the Auslander correspondence is the higher Auslander--Solberg correspondence\cite{IS18}. This is a correspondence between $n$-precluster tilting modules and Gorenstein algebras $\Gamma$ satisfying 
$$\id_\Gamma \Gamma \leq n+1 \leq \domdim \Gamma.$$ 
Such algebras are called $n$-minimal Auslander--Gorenstein algebras.

Precluster tilting can be defined as a weakening of $n$-cluster tilting, where instead of requiring ${}^{\perp_{n-1}}\M = \M = \M^{\perp_{n-1}}$ we only require ${}^{\perp_{n-1}}\M = \M^{\perp_{n-1}}$. This definition also works in the exact case and is discussed in \cref{sec:precluster-tilting}. However, this definition does not work when $n=1$.

In \cite{IS18} Iyama and Solberg define precluster tilting using the higher Auslander--Reiten translate, by requiring that $\tau_n(\M) \subseteq \M$ and $\tau_n^-(\M) \subseteq \M$. This is not available in the exact setting, so we replace it by the condition of $\M$ being relative $(n-1)$-cotilting (\cref{def:precluster_tilting}). This definition is motivated by the proof of \cite[Thm~4.4]{IS18}

A finite dimensional algebra is Gorenstein if and only if the subcategory of projectives is cotilting. This is therefore the restriction we put on $n$-minimal Auslander--Gorenstein categories. The analog of $\id_\Gamma \Gamma$ is simply the injective dimension of the projectives, while for the dominant dimension we use a similar construction as for Morita--Tachikawa categories (\cref{def:Auslander-Gorenstein}) .

For finite dimensional algebras the criteria for $\Gamma$ to be Gorenstein is redundant, since it can be proven that any algebra satisfying 
$$\id_\Gamma \Gamma \leq n+1 \leq \domdim \Gamma$$  
is Gorenstein. This eludes us to a weaker definition of $n$-minimal Auslander--Gorenstein categories, namely if we only require the projectives to be partial cotilting. We call this partial $n$-minimal Auslander--Gorenstein categories, and we similarly define partial $n$-precluster tilting subcategories to only be partial relative cotilting.

Our second main theorem is then two analogs of the Auslander--Solberg correspondence.

\begin{mainthm}[\cref{thm:partial_auslander_solberg_correspondence}, \cref{thm:Auslander-Solberg_correspondence}]
  There is a correspondence between equivalence classes of (partial) $n$-precluster tilting subcategories and equivalence classes of (partial) $n$-minimal Auslander--Gorenstein categories.
\end{mainthm}

In \cref{sec:comparison_to_classical_case} we verify that the correspondence actually is an analog of the Auslander--Solberg correspondence by comparing it to the classical case. Specifically, we show that a module $M$ over an Artin algebra, is $n$-precluster tilting if and only if $\add M$ is an $n$-precluster tilting subcategory. We show that $\mod\Gamma$ is $n$-minimal Auslander--Gorenstein if and only if $\Gamma$ is $n$-minimal Auslander--Gorenstein. And we show that the correspondence in \cref{sec:Auslander-Solber-correspondence} restricts to the classical one for Artin algebras, as described in \cite{IS18}.

In \cref{sec:projectives_and_injectives} we extend some propositions from \cite{henrard2020auslanders} about injective and projective objects to Morita--Tachikawa categories.

Lastly, in \cref{sec:examples}, we give a few simple examples of $n$-precluster tilting subcategories. 

\textbf{Acknowledgments} We would like to thank Øyvind Solberg for excellent guidance and input throughout this project. We would also like to thank Sondre Kvamme for help with finding examples of precluster tilting subcategories.

\section{Preliminaries}\label{sec:preliminaries}

\textbf{Convention} All subcategories are assumed to be full, and closed under isomorphism, direct sums and retractions.

\textbf{Note} This paper is based on \cite{henrard2020auslanders} and \cite{ebrahimi2021higher}, where it is assumed that all categories are essentially small. However, these assumptions are not necessary for the relevant theorems we use, so we make no such assumption.  

We begin this paper with a few short definitions and propositions about exact categories.

\begin{definition}
  If $\mathcal X \subseteq \E$ is a subcategory of an exact category we define subcategories 
  \begin{alignat*}{3}
    {}^\perp \mathcal X &:= \{ Y \in \E \mid \Hom(Y, X)&=0,&\; \forall X \in \mathcal X \}\\
    {}^{\perp_n} \mathcal X &:= \{ Y \in \E \mid \Ext^{i}(Y, X)&=0,&\; \forall X \in \mathcal X, 1 \leq i \leq n \}\\
    {}^{\perp_\infty} \mathcal X &:= \{ Y \in \E  \mid  \Ext^{i}(Y, X)&=0,&\; \forall X \in \mathcal X, \forall i>0\}\\
  \end{alignat*}
  Note that when $n=0$ the condition for ${}^{\perp_0} \mathcal X$ is void and so we have ${}^{\perp_0} \mathcal X = \E$. Dually we also define $\mathcal X^\perp$, $\mathcal X^{\perp_n}$ and $\mathcal X^{\perp_\infty}$. If $\mathcal X = \add X$ for an object $X$ then we may write $X^\perp$ for $\mathcal X^\perp$.
\end{definition}

\begin{definition}
  Let $\M \subseteq \E$ be a subcategory of an exact category and $X$ an object in $\E$. We define an \textit{$\M$-resolution} of $X$ to be an exact sequence 
  \begin{center}
    \begin{tikzcd}
      \cdots \ar[r] & M_2 \ar[r] & M_1 \ar[r] & M_0 \ar[r] & X \ar[r] & 0
    \end{tikzcd}
  \end{center}
  with $M_i \in \M$ such that each map is a right $\M$-approximation onto its image. I.e. the sequence remains exact when we apply $\Hom(M,-)$ for any $M$ in $\M$.

  Similarly an \textit{$\M$-presentation} of $X$ is an exact sequence 
  \begin{center}
    \begin{tikzcd}
      M_1 \ar[r] & M_0 \ar[r] & X \ar[r] & 0
    \end{tikzcd}
  \end{center}
  satisfying the same properties.

  We define \textit{$\M$-coresolution} and \textit{$\M$-copresentation} dually.
\end{definition}

\begin{definition}
  If $\M \subseteq \E$ is a subcategory of an exact category, we define the relative exact structure $F_\M$ by keeping those conflations that remain exact when we apply $\Hom(M, -)$ for any $M\in\M$. We define the structure $F^\M$ dually.
\end{definition}

Note that the above definition does indeed give us an exact structure by \cite[Prop~1.7]{DRSS99}. 

\begin{proposition}\cite[Prop~1.8]{AS93}\cite[III~Corollary~4.2]{Aus78}\label{prop:dual_of_relative_exact_structure}
  Let $\E$ be an exact category which has almost split sequences and let $\M$ be a subcategory. Then the two structures $F_\M$ and $F^{\tau M}$ coincide.
  \begin{proof}
    Let $\eta\colon 0 \to A \to B \to C \to 0$ be an $F_\M$-exact sequence. For the sake of contradiction, assume it is not $F^{\tau \M}$-exact. Then there is a map $f\colon A \to \tau M$ such that $f\cdot \eta \neq 0$. By the property of almost split sequences this means there is a map $g\colon M \to C$ such that $(f\cdot \eta) \cdot g$ is almost split. But $(f\cdot \eta) \cdot g = f\cdot (\eta\cdot g)$, and since $\eta$ is $F_\M$-exact $\eta\cdot g = 0$ for all $g$. This is a contradiction, and so the sequence is $F^{\tau \M}$-exact. The converse is proved similarly.
  \end{proof}
\end{proposition}

In this paper we consider cotilting in categories that don't necessarily have enough injectives. For this we need a more general definition of cotilting.

\begin{definition}[Cotilting] \cite[Def~7]{ZZ20}\label{def:cotilting}
  Let $\mathcal T \subseteq \E$ be a full additive subcategory of an exact category, closed under direct summands. We say that $\mathcal T$ is $n$-\emph{cotilting} iff 
  \begin{enumerate}
    \item $\mathcal T$ is functorinally finite,
    \item the injective dimension of $T$ is less than or equal to $n$ for every $T \in \mathcal T$,
    \item $\mathcal T$ is a cogenerator in ${}^{\perp_\infty} \mathcal T$. I.e. $\mathcal T \subseteq {}^{\perp_\infty} \mathcal T$ and for any $X\in{}^{\perp_\infty} \mathcal T$, there is a conflation $X \rightarrowtail T \twoheadrightarrow X'$ with $T \in \mathcal T$ and $X' \in {}^{\perp_\infty} \mathcal T$.
  \end{enumerate}

  Dually we say that $\mathcal T$ is \emph{$n$-tilting} if the projective dimension is bounded by $n$, it's functorially finite, and $\mathcal T$ is a generator in $\mathcal T^{\perp_\infty}$.
\end{definition}

\section{The category of admissibly presented modules}\label{sec:admissible_modules}

In this section we recall the category of admissibly presented modules introduced in \cite{henrard2020auslanders}, and some key properties that is used in the following sections.

\begin{definition}
  Let $\M \subseteq \E$ be a subcategory of an exact category. We use $\Mod(\M)$ to denote the category of additive contravariant functors from $\M$ to $\operatorname{Ab}$. We define the admissible module category $\mod_{\adm}(\M)$ to be the full subcategory of $\Mod(\M)$ consisting of functors $F$ with a projective presentation
  \begin{center}
    \begin{tikzcd}
      (-,M') \ar[r]{}{(-, f)} & (-,M) \ar[r] & F \ar[r] & 0
    \end{tikzcd}
  \end{center}
  where $f\colon M' \to M$ is admissible, and $(-,X):= \Hom_\E(-,X)|_\M$.
\end{definition}

\textbf{Remark} If $\M$ is not essentially small then $\Mod(\M)$ might not be locally small. However, $\mod(\M)$ and consequently $\mod_{\adm}(\M)$ are always locally small, so this does not create any set theoretic issues.

\begin{proposition}\cite[Prop~3.5]{henrard2020auslanders}
  The subcategory $\mod_{\adm}(\M) \subseteq \Mod(\M)$ is extension-closed, thus it inherits an exact structure.
\end{proposition}

We now prove a slight generalization of the Yoneda Lemma, which will be useful later.

\begin{lemma}[Yoneda Lemma]\label{lemma:yoneda}
  Let $\E$ be an exact category and let $\M \subseteq \E$ be a generating contravariantly finite subcategory. Then for any objects $X, Y \in \E$ we have that
  \begin{align*}
    \Hom_\E(X, Y) \cong \Hom_{\mod_{\adm}(\M)}((-,X), (-, Y)).
  \end{align*}
  \begin{proof}
    Take an $\M$-presentation of $X$ to be 
    \begin{center}
      \begin{tikzcd}
        M' \ar[r] & M \ar[r] & X \ar[r] & 0.
      \end{tikzcd}
    \end{center}
    This gives rise to a projective presentation of $(-,X)$:
    \begin{center}
      \begin{tikzcd}
        (-,M') \ar[r] & (-, M) \ar[r] & (-, X) \ar[r] & 0.
      \end{tikzcd}
    \end{center}
    Now if we apply $\Hom_{\mod_{\adm}(\M)}(-, (-, Y))$ we get 
    \begin{center}
      \begin{tikzcd}
        0 \ar[r] & ((-,X), (-, Y)) \ar[d, <-, dashed] \ar[r] & ((-, M), (-,Y)) \ar[d, <-, "\cong"] \ar[r] & ((-, M'), (-,Y)) \ar[d, <-, "\cong"]\\
        0 \ar[r] & (X, Y) \ar[r] & (M, Y) \ar[r] & (M', Y)
      \end{tikzcd}
    \end{center}
    where the two rightmost maps are isomorphisms by the Yoneda Lemma. Then the 5-lemma gives us that the dashed arrows is an isomorphism, which was our desired result.
  \end{proof}
\end{lemma}

We now recall some important subcategories of $\mod_{\adm}(\M)$ from \cite{henrard2020auslanders,ebrahimi2021higher} which we use when defining our correspondence in \cref{sec:Auslander-Solber-correspondence}.

\begin{definition}
  Let $\M \subseteq \E$ be a subcategory of an exact category. We define the following subcategories of $\mod_{\adm}(\M)$:
  \begin{enumerate}
    \item Let $\operatorname{eff}(\M)$ be the subcategory consisting of effaceable functors, i.e. those functors with a projective presentation 
    \begin{center}
      \begin{tikzcd}
        (-,M') \ar[r]{}{(-, f)} & (-,M) \ar[r] & F \ar[r] & 0
      \end{tikzcd}
    \end{center}
    where $f$ is a deflation.
    \item Let $\mathcal F(\M)$ be the subcategory of functors that admit a projective presentation 
    \begin{center}
      \begin{tikzcd}
        (-,M') \ar[r]{}{(-, f)} & (-,M) \ar[r] & F \ar[r] & 0
      \end{tikzcd}
    \end{center}
    where $M' \twoheadrightarrow \Image(f)$ is a right $\M$-approximation.
  \end{enumerate}
\end{definition}

\begin{proposition}
  Let $\M \subseteq \E$ be a subcategory of an exact category. If $\M$ is generating and contravariantly finite, then 
  \begin{enumerate}
    \item $\operatorname{eff}(\M)$ is equal to the full subcategory of $\Mod(\M)$ that consists of functors which fit into an exact sequence 
    \begin{center}
      \begin{tikzcd}
        (-,M) \ar[r]{}{(-, f)} & (-,C) \ar[r] & F \ar[r] & 0
      \end{tikzcd}
    \end{center}
    where $M$ is in $\M$, $C$ is in $\E$, and $f\colon M \to C$ is a deflation.
    \item $\mathcal F(\M)$ is equal to the full subcategory of $\Mod(\M)$ that consist of functors which fit into an exact sequence 
    \begin{center}
      \begin{tikzcd}
        (-,K) \ar[r]{}{(-, f)} & (-,M) \ar[r] & F \ar[r] & 0
      \end{tikzcd}
    \end{center}
    where $M$ is in $\M$, $K$ is in $\E$, and $f\colon K \to M$ is an inflation.
  \end{enumerate}
  \begin{proof}
    \begin{enumerate}
      \item[]
      \item Assume $F$ fits into an exact sequence of the desired form, and let the sequence be given by
      \begin{center}
        \begin{tikzcd}
          (-,M) \ar[r]{}{(-, f)} & (-,C) \ar[r] & F \ar[r] & 0.
        \end{tikzcd}
      \end{center}
      Let $M_C \to C$ be a right $\M$-approximation of $C$, and take a pullback along $f$.
      \begin{center}
        \begin{tikzcd}
          PB \arrow[dr, phantom, "\scalebox{1}{$\lrcorner$}" , very near start, color=black] \ar[d, twoheadrightarrow] \ar[r, twoheadrightarrow] & M_C \ar[d, twoheadrightarrow]\\
          M \ar[r, twoheadrightarrow, "f"] & C
        \end{tikzcd}
      \end{center}
      Applying the Yoneda embedding gives us the diagram
      \begin{center}
        \begin{tikzcd}
          (-,PB) \arrow[dr, phantom, "\scalebox{1}{$\lrcorner$}" , very near start, color=black] \ar[d] \ar[r] & (-,M_C) \ar[d, twoheadrightarrow] \ar[r] & F \ar[d, equal] \ar[r] & 0\\
          (-,M) \ar[r]{}{(-,f)} & (-,C) \ar[r] & F \ar[r] & 0
        \end{tikzcd}
      \end{center}
      Since the pullback along an epimorphism preserves cokernels we get that the top row is exact. Now, let $M_{PB} \twoheadrightarrow PB$ be an approximation of $PB$. Since $\M$ is generating we can choose this approximation to be a deflation (c.f. the dual of \cite[Prop~2.12(ii)]{Buhler10}), and we get that $F$ has a projective presentation coming from a deflation. Thus $F$ is in $\operatorname{eff}(\M)$.
      \item If $F$ has a presentation coming from a morphism $f\colon M' \to M$ where $M'$ is an approximation onto the image, we can choose $K$ to be the image of $f$. Conversely if $f$ comes from a sequence induced by an inflation we can make it into a projective presentation by taking an approximation of $K$.
    \end{enumerate}
  \end{proof}
\end{proposition}

We conclude this section by stating two propositions that were proved in \cite{henrard2020auslanders}.

\begin{proposition}\cite[Prop~3.6]{henrard2020auslanders}\label{prop:eff_F_torsion_pair}
  Let $\M \subseteq \E$ be a subcategory of an exact category. If $\M$ is generating and contravariantly finite. Then the pair $(\operatorname{eff}(\M), \mathcal F(\M))$ is a torsion pair in $\mod_{\adm}(\M)$ and $\operatorname{eff}(\M)$ satisfies axiom 
  {\bfseries A2} 
  of \cite[Def~6.1]{henrard2020localizations}. I.e. any morphism $F \to E$ in $\mod_{\adm}(\M)$ with $E\in \operatorname{eff}(\M)$ is admissible with image in $\operatorname{eff}(\M)$.
\end{proposition}

\begin{proposition}\cite[Prop~3.17]{henrard2020auslanders}\label{prop:eff_is_ort_P_F_is_cogenP}
  Let $\M \subseteq \E$ be a subcategory of an exact category. Let $\P$ denote the projective objects in $\mod_{\adm}(\M)$, and let $\eff$ be the left $\Hom$-orthogonal complement of $\P$. If $\M$ is generating-cogenerating and contravariantly finite, then
  \begin{enumerate}
    \item $\operatorname{eff}(\M) = \eff$, 
    \item $\mathcal F = \cogen \P$,
    \item The Yoneda embedding gives an equivalence from $\M$ to $\P$.
  \end{enumerate}
\end{proposition}

\section{Auslander--Solberg correspondence}\label{sec:Auslander-Solber-correspondence}

In this section we apply the theory of localization at two-sided admissibly percolating subcategories as described in \cite{henrard2020localizations} to construct a correspondence between Morita--Tachikawa categories and contravariantly finite generating--cogenerating subcategories. This restricts to a correspondence between (partial) $n$-minimal Auslander--Gorenstein categories and (partial) $n$-precluster tilting subcategories. We use the term \emph{weak equivalence} to refer to a map which becomes an isomorphism under localization, and we write $X \xrightarrow{\sim} Y$ for a weak equivalence from $X$ to $Y$.

\subsection{The bijection}

In \cite{henrard2020auslanders} a correspondence between \emph{Auslander exact categories} and exact categories is constructed. Here we generalize this correspondence slightly to what we have chosen to call \emph{Morita--Tachikawa categories}.

\begin{definition}[\pac]\label{def:pac}
  Let $\E$ be an exact category and let $\P$ the projective objects of $\E$. We write $\eff$ for the $\Hom$-orthogonal complement $$\eff := \{X \in \E | \Hom_\E(X, P)=0, \forall P \in \P\}.$$ We say that $\E$ is a \emph{\pac} iff
  \begin{enumerate}
    \item $\E$ has enough projectives.
    \item $(\eff, \cogen \P)$ forms a torsion pair.
    \item \label{item:axiom_A2} $\eff$ satisfies axiom \textbf{A2} of \cite[Def~6.1]{henrard2020localizations}, i.e. any morphism $X \to E$ with $E\in \eff$ is admissible with image in $\eff$.
    \item $\Ext^{1}(\eff, \P)=0$.
  \end{enumerate}
  Note, condition 2 and 3 together gives us that $\eff$ is two-sided admissibly percolating \cite[Prop~2.22(1)]{henrard2020auslanders}.
\end{definition}

Throughout this section we consider the localization $\E/\eff$. The following lemma ensures that this is a well defined category even when we don't assume $\E$ to be essentially small.

\begin{lemma}\label{lem:equiv_on_syzygy}
  Let $\E$ be a \pac, with projective objects $\P$. Consider the localization functor $\E \to \E/\eff$. Then for every $X \in \E$ there is a weak equivalence $X \xrightarrow{\sim} F$ where $F$ is a second syzygy. Further, whenever $F$ is a second syzygy and $X$ is arbitrary, we have $$\Hom_{\E}(X, F) \cong \Hom_{\E / \eff}(X, F).$$ In particular the localization functor induces an equivalence of (non-exact) categories between $\Omega^2 \E$ and $\E/\eff$. Further, the inverse $\E/\eff \to \Omega^2 \E \subseteq \E$ is left exact, and gives a right adjoint to the localization functor.
  \begin{proof}
    Since $(\eff, \cogen \P)$ is a torsion pair we have a conflation $tX \rightarrowtail X \twoheadrightarrow fX$ with $tX \in \eff$ and $fX \in \cogen \P$. This shows that $X$ is weakly equivalent to a first syzygy, namely $fX$. The argument in \cite[Prop~4.4(3)]{henrard2020auslanders} shows that every first syzygy is weakly equivalent to a second syzygy, and thus $X$ is weakly equivalent to a second syzygy.

    Next we show that $\Hom_{\E}(X, F) \cong \Hom_{\E / \eff}(X, F)$ whenever $F$ is a second syzygy. The argument is similar to \cite[Prop~4.5]{henrard2020auslanders}. The localization functor gives a map $\Hom_{\E}(X, F) \to \Hom_{\E / \eff}(X, F)$. We wish to show that this is a bijection. Let $f\colon X \to F$ be in the kernel. Then by \cite[Prop~2.22(1)]{henrard2020auslanders} and \cite[Prop~4.10]{henrard2020localizations} there is a weak equivalence $g \colon F \xrightarrow{\sim} Y$ such that the composition $gf$ is 0. The kernel of $g$ is an element of $\eff$ mapping to $F$, hence it must be 0. Thus $g$ is a monomorphism and it follows that $f=0$, which shows injectivity.

    We now show surjectivity. Any map in $\Hom_{\E / \eff}(X, F)$ is of the form $X \to Y \xleftarrow{\sim} F$ \cite[Prop~4.6, Construction~2.8(3)]{henrard2020localizations}. We want to show that the weak equivalence is a split monomorphism. Just like before it is clear that it is mono, so we have a conflation $F \rightarrowtail Y \twoheadrightarrow tY$, with $tY \in \eff$. Since $F$ is a second syzygy, we also have a conflation $F \rightarrowtail P \twoheadrightarrow F'$ with $P \in \P$ and $F' \in \cogen \P$. Then we get an exact sequence 
    \begin{center}
      \begin{tikzcd}
        0=\Hom(tY, F') \ar[r] & \Ext^1(tY, F) \ar[r] & \Ext^{1}(tY, P)=0.
      \end{tikzcd}
    \end{center}
    Thus $\Ext^1(tY, F)=0$ and the map is split. Then we get a preimage of $X \to Y \xleftarrow{\sim} F$ by composing $X\to Y$ with the splitting. Thus the map is surjective and so 
    $$\Hom_{\E}(X, F) \cong \Hom_{\E / \eff}(X, F).$$

    Lastly, since when we restrict the localization to $\Omega^2 \E$ we get a dense fully faithful functor, this means that it is an equivalence between $\Omega^2\E$ and $\E/\eff$. The formula 
    $$\Hom_{\E}(X, F) \cong \Hom_{\E / \eff}(X, F)$$
    says exactly that the inverse is right adjoint to the localization. Right adjoints preserve limits, so they are in particular left exact.
  \end{proof}
\end{lemma}

With this lemma we now have the necessary tools to construct our desired bijection.

\begin{theorem}[The exact Morita--Tachikawa correspondence]\label{thm:bijection_partial_auslander}
  There is a bijection between 
  \begin{enumerate}
    \item equivalence classes of Morita--Tachikawa categories,
    \item equivalence classes of pairs consisting of an exact category and a contravariantly finite generating-cogenerating subcategory,
  \end{enumerate}
  where equivalence means equivalence of exact categories, and in the latter case we require the equivalence to restrict to an equivalence between the subcategories.
  
  The bijection is given by mapping a \pac{} $\E$ with projectives $\P$ to the image of $\P$ in the localization $\E/\eff$, and mapping an exact category $\E'$ with subcategory $\M$ to $\mod_{\adm}(\M)$.
  \begin{proof}
    Let $\E$ be a \pac, and let $\M$ be the image of the projectives under the localization $\E \to \E/\eff$. Then because $\E$ has enough projectives and projectives are closed under retractions we see that $\M$ is generating and closed under retractions. By \cref{lem:equiv_on_syzygy} every object in $\E/\eff$ is the image of a second syzygy. Since second syzygies embeds into projectives $\M$ is cogenerating. Lastly, since \cref{lem:equiv_on_syzygy} gives us an equivalence, we have that an $\M$-approximation corresponds to a $\P$-approximation of a second syzygy. Since $\P$ is contravariantly finite, so is $\M$.

    For the other direction we see that $\mod_{\adm}(\M)$ has enough projectives, simply by how it is defined.

    That $(\eff, \cogen\P)$ forms a torsion pair, and that $\eff$ satisfies \textbf{A2} is the statement of \cref{prop:eff_F_torsion_pair} and \cref{prop:eff_is_ort_P_F_is_cogenP}.

    We have seen in \cref{prop:eff_is_ort_P_F_is_cogenP} that if $C$ is in $\eff$, then it has a presentation 
    \begin{center}
      \begin{tikzcd}
        0 \ar[r] & (-, K) \ar[r] & (-,M_1) \ar[r] & (-,M_0) \ar[r] & C \ar[r] & 0
      \end{tikzcd}
    \end{center}
    where $M_1 \to M_0$ is a deflation. Let $M_2$ be an approximation of $K$. Then 
    \begin{center}
      \begin{tikzcd}
        M_2 \ar[r] & M_1 \ar[r] & M_0 \ar[r] & 0
      \end{tikzcd}
    \end{center}
    is an exact sequence. Applying $(-,M)$ we get an exact sequence 
    \begin{center}
      \begin{tikzcd}
        0 \ar[r] & (M_0, M) \ar[r] & (M_1, M) \ar[r] & (M_2, M).
      \end{tikzcd}
    \end{center}
    The homology at $M_1$ is 0 and equal to $\Ext^1(C, (-,M))$, thus $\Ext^{1}(\eff, \P)=0$.

    To see that the two mappings are inverses of each other let $\P$ be the projectives in $\E$ and let $\M$ be its image under the localization $\E \to \E/\eff$. By \cref{lem:equiv_on_syzygy} the localization gives an equivalence $\P \cong \M$. By \cite[Thm~2.16, Prop~2.22(1)]{henrard2020auslanders} this equivalence preserves and reflects admissible morphisms, thus $\mod_{\adm}(\P) \cong \mod_{\adm}(\M)$. Lastly, by \cite[Prop~3.5]{ebrahimi2021higher} we have $\E \cong \mod_{\adm}(\P)$.

    By \cref{prop:eff_is_ort_P_F_is_cogenP} $\M$ is equivalent to the category of projectives in $\mod_{\adm}(\M)$ and by \cite[Prop~3.15]{ebrahimi2021higher} $\mod_{\adm}(\M)/\eff \cong \E'$, so the result follows.
  \end{proof}
\end{theorem}

\subsection{The partial Auslander--Solberg correspondence}

In this section we define the notion of (partial) $n$-minimal Auslander--Gorenstein categories and (partial) $n$-precluster tilting subcategories. In \cref{thm:partial_auslander_solberg_correspondence} and \cref{thm:Auslander-Solberg_correspondence} we show that the bijection constructed in \cref{thm:bijection_partial_auslander} restricts to a bijection between these.

\begin{definition}[$n$-minimal Auslander--Gorenstein]\label{def:Auslander-Gorenstein}
  Let $\E$ be a \pac{} and use the same notation as in \cref{def:pac}. We say that $\E$ is \emph{partial $n$-minimal Auslander--Gorenstein} if it satisfies
  \begin{enumerate}
    \item \label{item:proj_covar_fin} $\P$ is covariantly finite with admissible left-approximations.
    \item \label{item:domdim} $\Ext^{<n+1}(\eff, \P)=0$.
    \item \label{item:proj_fin_injdim} The injective dimension of any projective is less than or equal to $n+1$, i.e. $\Ext^{>n+1}(-, \P)=0$.
  \end{enumerate}
  We say that $\E$ is \emph{$n$-minimal Auslander--Gorenstein} if it in addition satisfies
  \begin{enumerate}[resume]
    \item \label{item:proj_cotilt} The subcategory $\P$ is cotilting.
  \end{enumerate}
\end{definition}

\begin{definition}[$n$-precluster tilting]\label{def:precluster_tilting}
  Let $\M$ be a subcategory of an exact category $\E'$. We say that $\M$ is \emph{partial $n$-precluster tilting} if it satisfies
  \begin{enumerate}[label=(\alph*)]
    \item \label{item:gen-cogen} $\M$ is a generating-cogenerator subcategory.
    \item \label{item:contrafin} $\M$ is contravariantly finite.
    \item \label{item:covarfin} $\M$ is covariantly finite.
    \item \label{item:n-rigid} $\M$ is $n$-rigid, i.e. $\Ext^{<n}(\M, \M) = 0$. 
    \item \label{item:relative_injdim} The relative injective dimension $\id_{F_{\M}} M$ is strictly less than $n$ for any $M\in\M$.
  \end{enumerate}
  We say that $\M$ is \emph{$n$-precluster tilting} if it in addition satisfies
  \begin{enumerate}[resume, label=(\alph*)]
    \item \label{item:relative_cotilting}$\M$ is $F_\M$-cotilting.
  \end{enumerate}
\end{definition}

We shall now see how $n$-minimal Auslander--Gorenstein categories and $n$-precluster tilting subcategories interact with the bijection defined in \cref{thm:bijection_partial_auslander}.

\begin{theorem}\label{thm:partial_n-AG_maps_to_almsot_n-precluster}
  If $\E$ is partial $n$-minimal Auslander--Gorenstein, then $\M:=\P \subseteq \E/\eff$ is partial $n$-precluster tilting.
  \begin{proof}
    \begin{itemize}[leftmargin=*]
    \item[]
    \item[] Since $\E$ is a \pac{} we have that $\M$ is generating-cogenerating and contravariantly finite, hence condition \ref{item:gen-cogen} and \ref{item:contrafin} of \cref{def:precluster_tilting} hold. We proceed by proving conditions \ref{item:covarfin}, \ref{item:n-rigid}, and \ref{item:relative_injdim}.
    \item[\ref{item:covarfin}]
    By \cref{lem:equiv_on_syzygy} we have that $\E/\eff$ is equivalent to $\Omega^2\E$. Since $\P$ is covariantly finite in $\Omega^2\E$, this means that $\M$ is covariantly finite in $\E/\eff$.

    \item[\ref{item:n-rigid}]
    The fact that $\Ext^{<n}(\M, \M) = 0$ follows from a similar argument to \cite[Thm~4.3]{ebrahimi2021higher}. 
    
    When $n=1$ there is nothing to show, so assume $n>1$. We want to show that $\Ext^k_{\E/\eff}(\M, \M)=0$ for every $k$ with $1 \leq k < n$. We do this by induction on $k$. If $k=1$ an element of $\Ext^1_{\E/\eff}(\M, \M)$ is a sequence of the form 
    \begin{center}
      \begin{tikzcd}
        0 \ar[r] & M' \ar[r] & X_1 \ar[r] & M \ar[r] & 0
      \end{tikzcd}
    \end{center}
    If we denote the localization functor $\E \to \E/\eff$ by $H$, then we know $M = H(P^0)$ and $M'=H(P)$ for projectives $P^0$ and $P$. Further by \cref{lem:equiv_on_syzygy} the map $X_1 \to M$ comes from an exact sequence $X_1 \to P^0 \to C \to 0$ with $C$ in $\eff$. By choosing $P^1$ to be a projective mapping onto $X_1$ we can construct a projective resolution of $C$
    \begin{center}
      \begin{tikzcd}
        P^\bullet \colon \quad P^{n+1} \ar[r] & \cdots \ar[r] & P^1 \ar[r] & P^0 \ar[r] & C \ar[r] & 0.
      \end{tikzcd}
    \end{center}
    Since $\Ext^{<n+1}(\eff, \P)=0$ this sequence remains exact if we apply $\Hom_\E(-, P)$. Further, since past $P^1$ all images are second syzygies, \cref{lem:equiv_on_syzygy} gives us that the sequence we get by applying the localization,
    \begin{center}
      \begin{tikzcd}
        H(P^{n+1}) \ar[r] & H(P^n) \ar[r] & \cdots \ar[r] & H(P^2) \ar[r] & H(P^1),
      \end{tikzcd}
    \end{center}
    remains exact when we apply $\Hom_{\E/\eff}(-, H(P))=\Hom_{\E/\eff}(-, M')$. Now the mapping $P^1 \to X_1$ induces the diagram
    \begin{center}
      \begin{tikzcd}
        H(P^3) \ar[d] \ar[r, "d_P^3"] & H(P^2) \ar[r, "d_P^2"] \ar[d, "f"] & H(P^1) \ar[r] \ar[dl, dashed] \ar[d, "\pi"] & H(P^0) \ar[d, equal] \ar[r] & 0\\
        0 \ar[r] & M' \ar[r, "\iota"] & X_1 \arrow[ul, phantom, "\scalebox{1}{$\ulcorner$}" , very near start, color=black] \ar[r] & M \ar[r] & 0.
      \end{tikzcd}
    \end{center}
    A diagram chase shows that the left square is pushout as indicated\cite[Prop~2.12]{Buhler10}. Since $\iota f d_P^3=\pi d_P^2 d_P^3 =0$ and $\iota$ is a monomorphism, we have $fd_P^3=0$. Now since $\Hom_{\E/\eff}(H(P^\bullet), M')$ is exact we have that $f$ factors through $H(P^1)$ along the dotted arrow. By the pushout property the bottom sequence splits, thus $$\Ext^1_{\E/\eff}(\M, \M) = 0.$$

    Now let $k$ be an integer $1<k<n$ and assume by induction that $\Ext^{<k}_{\E/\eff}(\M, \M) = 0$. We consider an extension in $\Ext^k_{\E/\eff}(\M, \M)$ of the form 
    \begin{center}
      \begin{tikzcd}
        0 \ar[r] & M' \ar[r] & X_k \ar[r] & \cdots \ar[r] & X_1 \ar[r] & M \ar[r] & 0
      \end{tikzcd}
    \end{center}
    Just like before we construct the complex $P^\bullet$. Taking the pullback of $H(P^1) \to X_1$ along $X_2 \to X_1$ yields a Yoneda-equivalent extension
    \begin{center}
      \begin{tikzcd}[column sep=18pt]
        0 \ar[r] & M' \ar[d, equal] \ar[r] & X_k \ar[d, equal] \ar[r] & \cdots \ar[r] & X_3 \ar[d, equal] \ar[r] & PB \ar[d] \ar[r]\arrow[dr, phantom, "\scalebox{1}{$\lrcorner$}" , very near start, color=black] & H(P^1) \ar[d]\ar[r] & H(P^0)\ar[d, equal] \ar[r] & 0\\
        0 \ar[r] & M' \ar[r] & X_k \ar[r] & \cdots \ar[r] & X_3 \ar[r] & X_2 \ar[r] & X_1 \ar[r] & M \ar[r] & 0.
      \end{tikzcd}
    \end{center}
    Let $K_i$ denote the kernel of $H(P^i) \to H(P^{i-1})$. Then the above extension factors as
    \begin{center}
      \begin{tikzcd}[column sep=18pt]
        0 \ar[r] & M' \ar[r] & X_k \ar[r] & \cdots \ar[r] & PB \ar[rr] \ar[dr, twoheadrightarrow] && H(P^1) \ar[r] & H(P^0) \ar[r] & 0.\\
        &&&&& K_1 \ar[ur, hookrightarrow]
      \end{tikzcd}
    \end{center}
     
    Thus for it to be trivial it is enough to show that $\Ext^{k-1}_{\E/\eff}(K_1, M')=0$. 
    
    Consider the short exact sequence $K_i \rightarrowtail H(P^i) \twoheadrightarrow K_{i-1}$, and apply the long exact sequence in $\Ext(-,M')$. Since $\Ext^{<k}_{\E/\eff}(\M, \M) = 0$ this gives us
    $$\Ext^{k-1}_{\E/\eff}(K_1, M') \; = \; \Ext^{k-2}_{\E/\eff}(K_2, M') \; =\;\cdots\; =\;\Ext^{1}_{\E/\eff}(K_{k-1}, M').$$

    Now since $\Hom_{\E/\eff}(H(P^\bullet), M')$ is exact the map $K_k \to H(P^k)$ is a left $\M$-approximation. Then we have another long exact sequence in $\Ext$
    \begin{center}
      \begin{tikzcd}[column sep=15pt]
        \Hom(H(P^{k}), M') \ar[r] & \Hom(K_{k}, M') \ar[r, "0"] & \Ext^{1}(K_{k-1}, M') \ar[r] & \Ext^{1}(H(P^{k}), M').
      \end{tikzcd}
    \end{center}
    Since $\Ext^{1}(H(P^{k}), M')=0$ this shows that $\Ext^{1}_{\E/\eff}(K_{k-1}, M')=0$, and we conclude that $\Ext^{<n}(\M, \M) = 0$. 

    \item[\ref{item:relative_injdim}]
    Lastly for an arbitrary $X$ we choose $F$ to be a second syzygy weakly equivalent to $X$. Then an $\M$-resolution of $X$ corresponds to a projective resolution of $F$. Let $M$ be an object in $\M$ and let $P$ be the corresponding projective in $\E$. Applying $\Hom_{\E}(-, P)$ to the resolution of $F$ we get a complex
    \begin{center}
      \begin{tikzcd}
        0 \ar[r] & \Hom_{\E}(P^0, P) \ar[r] & \Hom_{\E}(P^1, P) \ar[r] & \cdots
      \end{tikzcd}
    \end{center}
    whose homology is equal to $\Ext_\E^\bullet(F, P)$. By \cref{lem:equiv_on_syzygy} this complex is isomorphic to 
    \begin{center}
      \begin{tikzcd}
        0 \ar[r] & \Hom_{\E/\eff}(P^0, M) \ar[r] & \Hom_{\E/\eff}(P^1, M) \ar[r] & \cdots
      \end{tikzcd}
    \end{center}
    whose homology is $\Ext_{F_{\M}}^\bullet(X, M)$. Thus $\Ext_{\E}^\bullet(F, P) = \Ext_{F_{\M}}^\bullet(X, M)$. Now since $F$ is a second syzygy and $\id P \leq n+1$ we have $\Ext_{F_{\M}}^{>n-1}(X, M) = \Ext_{\E}^{>n-1}(F, P)=0$. Thus $\id_{F_{\M}} M \leq n-1$.
  \end{itemize}
  \end{proof}
\end{theorem}

We have now showed that the bijection of \cref{thm:bijection_partial_auslander} maps partial Auslander--Gorenstein categories to partial precluster tilting subcategories. Next we show the other direction.

\begin{theorem}\label{thm:almsot_n-precluster_maps_to_partial_n-AG}
  If $\M \subseteq \E$ is partial $n$-precluster tilting, then $\mod_{\adm} \M$ is partial $n$-minimal Auslander--Gorenstein.
  \begin{proof}
    \begin{itemize}[leftmargin=*]
      \item[]


    \item[(\ref{item:proj_covar_fin})]
    That the projectives in $\mod_{\adm}(\M)$ are admissibly covariantly finite is proven in \cite[Prop~3.19(2)]{ebrahimi2021higher}. 


    \item[(\ref{item:domdim})]
    Let $\P$ denote the projective objects in $\mod_{\adm}(\M)$. A similar argument to \cite[Prop~3.17(4)]{ebrahimi2021higher} gives that $\Ext^{<n+1}(\eff, \P)=0$. Let $X$ be in $\eff \subseteq \mod_{\adm}(\M)$. Then \cref{prop:eff_is_ort_P_F_is_cogenP} gives us that $X$ has a projective presentation coming from a deflation $M_1 \twoheadrightarrow M_0$. Since $\M$ is generating contravariantly finite, we can continue this to an exact sequence by taking right $\M$-approximations:
    \begin{center}
      \begin{tikzcd}
        M_{n+1}\ar[r] & M_n \ar[r] & \cdots \ar[r] & M_1 \ar[r] & M_0 \ar[r] & 0.
      \end{tikzcd}
    \end{center}
    Let $K_i$ be the kernel of $M_i \to M_{i-1}$. Since $\M$ is $n$-rigid we have $\Ext^1(K_{i}, M) = \Ext^{i}(K_1, M)$ for all $i < n$. Further, from the conflation $K_1 \rightarrowtail M_1 \twoheadrightarrow M_0$ we get that $\Ext^{i}(K_1, M)=0$ for all $i < n-1$. This means that for all $i<n$ the inclusion $K_i \to M_i$ is a left approximation. Thus the sequence remains exact when we apply $\Hom(-, M)$ and we have $\Ext^{<n+1}(X, (-,M))=0$.

    \item[(\ref{item:proj_fin_injdim})]
    Let $A$ be an object in $\E$. Since $\M$ is contravariantly finite and generating we can take an $\M$-resolution 
    \begin{center}
      \begin{tikzcd}
        \cdots \ar[r] & M_2 \ar[r] & M_1 \ar[r] & M_0 \ar[r] \ar[d, twoheadrightarrow] & 0\\
               & &    & A
      \end{tikzcd}
    \end{center}
    This induces a projective resolution of $\Hom(-, A)|_\M$ in $\mod_{\adm}(\M)$. For brevity we write $(-, A)$ instead of $\Hom(-, A)|_\M$.
    \begin{center}
      \begin{tikzcd}
        \cdots \ar[r] & (-,M_2) \ar[r] &(-,M_1) \ar[r] & (-,M_0) \ar[r] \ar[d, twoheadrightarrow] & 0\\
        &&& (-,A)
      \end{tikzcd}
    \end{center}
    To compute $\Ext^\bullet((-,A), (-, B))$ we apply $\Hom_{\mod_{\adm}(\M)}(-, (-, B))$. By the Yoneda lemma this gives us 
    \begin{center}
      \begin{tikzcd}
        0 \ar[r] & (M_0, B) \ar[r] & (M_1, B) \ar[r] & (M_2, B) \ar[r] & \cdots
      \end{tikzcd}
    \end{center}
    whose homology is equal to $\Ext_{F_{\M}}^\bullet(A, B)$. Thus we have 
    $$\Ext^\bullet((-,A), (-, B)) = \Ext_{F_{\M}}^\bullet(A, B).$$ 
    Now if $C$ is an arbitrary object in $\mod_{\adm} \M$, then it has a projective presentation 
    \begin{center}
      \begin{tikzcd}
        0 \ar[r] & (-, K) \ar[r] & (-, M'') \ar[r] & (-,M') \ar[r] & C \ar[r] & 0.
      \end{tikzcd}
    \end{center}
    Let $M$ be an object in $\M$. Because $\id_{F_{\M}} M \leq n-1$, we have $$\Ext^{>n+1}(C, (-, M)) = \Ext^{>n-1}((-,K), (-, M)) = \Ext^{>n-1}_{F_M}(K, M)=0.$$ Which concludes the proof.
  \end{itemize}
  \end{proof}
\end{theorem}

We now combine the two previous theorem into one main result.

\begin{theorem}[The exact partial Auslander--Solberg correspondence]\label{thm:partial_auslander_solberg_correspondence}
  The bijection given in \cref{thm:bijection_partial_auslander} restricts to a bijection between partial $n$-minimal Auslander--Gorenstein exact categories and partial $n$-precluster tilting subcategories.
  \begin{proof}
    \cref{thm:partial_n-AG_maps_to_almsot_n-precluster} and \cref{thm:almsot_n-precluster_maps_to_partial_n-AG} gives us the maps between the two classes and \cref{thm:bijection_partial_auslander} gives us that these maps are mutual inverses.
  \end{proof}
\end{theorem}

\subsection{The Auslander--Solberg correspondence}

In this section we strengthen \cref{thm:partial_auslander_solberg_correspondence} by restricting to $n$-minimal Auslander--Gorenstein categories and $n$-precluster tilting subcategories.

\begin{proposition}\label{prop:AG_implies_precluster}
  If $\E$ is $n$-minimal Auslander--Gorenstein, then $\M:=\P \subseteq \E/\eff$ is an $n$-precluster tilting subcategory.
  \begin{proof}
    We have already seen in \cref{thm:partial_n-AG_maps_to_almsot_n-precluster} that $\M$ is partial $n$-precluster tilting, so we just need to show that $\M$ is $F_\M$-cotilting.

    Let $X$ satisfy $\Ext^{>0}_{F_\M}(X, \M)=0$. We may assume $X$ is a second syzygy. Then 
    $$\Ext^{>0}_{F_\M}(X, \M)=\Ext^{>0}_\E(X, \P)=0.$$ 
    So $X$ is in ${}^{\perp_\infty}\P$ and, since $\P$ is cotilting, we have an exact sequence
    \begin{center}
      \begin{tikzcd}
        0 \ar[r] & X \ar[r] & P \ar[r] & C \ar[r] & 0
      \end{tikzcd}
    \end{center}
    with $C$ also in ${}^{\perp_\infty}\P$. Then $C$ is also a second syzygy, so by \cref{lem:equiv_on_syzygy} $P \to C$ induces an $\M$-approximation and 
    $$\Ext^{>0}_{F_\M}(C, M)=\Ext^{>0}(C, \P)=0.$$ 
    Thus $\M$ is $F_\M$-cotilting, which is what we wanted to prove.
  \end{proof}
\end{proposition}

\begin{proposition}\label{prop:precluster_implies_AG}
  If $\M \subseteq \E'$ is an $n$-precluster tilting subcategory, then $\mod_{\adm}(\M)$ is $n$-minimal Auslander--Gorenstein.
  \begin{proof}
    In \cref{thm:almsot_n-precluster_maps_to_partial_n-AG} we saw that $\mod_{\adm}(\M)$ is partial $n$-minimal Auslander--Gorenstein, so the only thing we need to show is that the projectives $\P \subseteq \mod_{\adm}(\M)$ forms a cotilting subcategory. To see this, let $X$ be an object such that $\Ext^{>0}(X, \P)=0$. Since $(\eff, \cogen\P)$ is a torsion pair we have conflation $tX \rightarrowtail X \twoheadrightarrow fX$ with $tX \in \eff$ and $fX\in\cogen\P$. Let $P$ be a projective and take the long exact sequence in $\Ext(-,P)$.
    \begin{center}
      \begin{tikzcd}
        \cdots \ar[r] & \Hom(tX, P) \ar[r] & \Ext^1(fX, P) \ar[r] & \Ext^1(X, P) \ar[dll, out=-15, in=165] \\ 
        &\Ext^1(tX, P) \ar[r] & \Ext^2(fX, P) \ar[r] & \Ext^2(X, P) \ar[r] & \cdots
      \end{tikzcd}
    \end{center}
    Since $\Hom(tX, P)=0$ and $\Ext^{>0}(X, P)=0$ this gives us that $\Ext^1(fX, P)=0$ and $\Ext^{i}(tX, P) = \Ext^{i+1}(fX, P)$ for all $i>0$. Since $\Ext^{<n+1}(\eff, \P)=0$ we have that $\Ext^{i}(tX, P)=0$ for all $i<n+1$. Since $\Ext^{>n+1}(-,\P)$ we have that $\Ext^{i}(tX, P) = \Ext^{i+1}(fX, P) = 0$ for all $i \geq n+1$. Thus $tX$ satisfies $\Ext^{>0}(tX, \P)=0$. 

    Remember that we want to show that $X$ embeds into a projective. To do this we need to show that $tX=0$. Let $C$ be an arbitrary object in $\eff$ that satisfies $\Ext^{>0}(C, \P)=0$. We now show that $C$ must be $0$.

    Take a projective presentation of $C$ to be 
    \begin{center}
      \begin{tikzcd}
        0 \ar[r] & \Omega^2 C \ar[r] & P_1 \ar[r] & P_0 \ar[r] & C \ar[r] & 0
      \end{tikzcd}
    \end{center}
    We know that $\Omega^2 C$ is isomorphic to $(-,K)$ for some $K$ in $\E'$. By the proof in \cref{thm:almsot_n-precluster_maps_to_partial_n-AG} we know that $\Ext_{F_{\M}}^{>0}(K, \M) = \Ext^{>0}(\Omega^2 C, \P) = \Ext^{>2}(C, \P)=0$. Then since $\M$ is cotilting we get an $F_\M$-exact sequence 
    \begin{center}
      \begin{tikzcd}
        0 \ar[r] & K \ar[r] & M_1 \ar[r] & M_2 \ar[r] & M_3 \ar[r] & \cdots
      \end{tikzcd}
    \end{center}
    Since this sequence is $F_\M$-exact, it remains exact when we apply the Yoneda embedding. Compare this sequence to our original presentation
    \begin{center}
      \begin{tikzcd}
        0 \ar[r] & \Omega^2 C \ar[d, equal] \ar[r] & P_1 \ar[r] \ar[d, dashed] & P_0 \ar[r] \ar[d, dashed] & C \ar[r] \ar[d, dashed, "0"] & 0\\
        0 \ar[r] & (-,K) \ar[r] & (-,M_1) \ar[r] & (-,M_2) \ar[r] & (-,M_3) \ar[r] & \cdots
      \end{tikzcd}
    \end{center}
    Since $\Ext^{>0}(C, \P)=0$, the map $\Omega^2 C \to P_1$ is a left $\P$-approximation. Thus we get a map $P_1 \to (-,M_1)$ making the left square commute. Similarly we get maps along the dashed arrows, making the diagram commute. Now, since $C$ is in $\eff$ the rightmost map must be 0. This means that $P_0 \to (-,M_2)$ factors through the kernel of $(-,M_2) \to (-,M_3)$, and so we get a homotopy
    \begin{center}
      \begin{tikzcd}
        0 \ar[r] & \Omega^2 C \ar[d, equal] \ar[r, "d"] & P_1 \ar[dl, dashed, swap, "h"] \ar[r] \ar[d] & P_0 \ar[r] \ar[dl, dashed] \ar[d] & C \ar[r] \ar[d, "0"] & 0\\
        0 \ar[r] & (-,K) \ar[r, "\iota"] & (-,M_1) \ar[r] & (-,M_2) \ar[r] & (-,M_3) \ar[r] & \cdots
      \end{tikzcd}
    \end{center}
    Then we have $\iota h d = \iota$ and since $\iota$ is mono this means that $d$ splits, so $\Omega C$ is projective. Since $\Ext^1(C, \P)=0$ this means that the resolution of $C$ splits, so $C=0$.

    In conclusion, we have shown that if $X$ is in ${}^{\perp_\infty}\P$, then $tX = 0$, which means that $X$ embeds into a projective. That is, there is a conflation $X \hookrightarrow P \twoheadrightarrow X'$. Since $\P$ is covariantly finite, we may choose $X \hookrightarrow P$ to be a left $\P$-approximation, and then it follows that $X'$ is also in ${}^{\perp_\infty}\P$. Thus, $\P$ is cotilting, and $\mod_{\adm}(\M)$ is $n$-minimal Auslander--Gorenstein.
  \end{proof}
\end{proposition}

Combining these two propositions gives us the main theorem of this section.

\begin{theorem}[The exact Auslander--Solberg correspondence]\label{thm:Auslander-Solberg_correspondence}
  The bijection given in \cref{thm:bijection_partial_auslander} restricts to a bijection between $n$-minimal Auslander--Gorenstein exact categories and $n$-precluster subcategories.
  \begin{proof}
    This is the combined statement of \cref{prop:AG_implies_precluster} and \cref{prop:precluster_implies_AG}.
  \end{proof}
\end{theorem}

We conclude this section by noting that $n$-Auslander exact categories are special cases of $n$-minimal Auslander--Gorenstein categories and that $n$-cluster tilting subcategories are special cases of $n$-precluster tilting subcategories.  

\begin{proposition}
  If $\E$ is $n$-Auslander then $\E$ is $n$-minimal Auslander--Gorenstein.
  \begin{proof}
    We only need to show that the projectives have injective dimension less than or equal to $n+1$, and that they are cotilting. Since $\E$ has finite global dimension less than or equal to $n+1$, the first assertion is immediate. Since $\E$ has finite global dimension and enough projectives, we can for any object $C$ find a $d$ such that $\Ext^d(C, \P) \neq 0$ by letting $d$ be the projective dimension of $C$. Thus ${}^{\perp_\infty} \P = \P$, so $\P$ is cotilting.
  \end{proof}
\end{proposition}
\begin{corollary}
  If $\M$ is $n$-cluster tilting, then $\M$ is $n$-precluster tilting.
  \begin{proof}
    By the correspondence in \cite{ebrahimi2021higher} we get that $\mod_{\adm}(\M)$ is $n$-Auslander, thus also $n$-minimal Auslander--Gorenstein. Then by \cref{thm:Auslander-Solberg_correspondence} $\M$ is $n$-precluster tilting.
  \end{proof}
\end{corollary}

\section{\texorpdfstring{\boldmath $n$}{n}-precluster tilting subcategories}\label{sec:precluster-tilting}

In the previous section we defined $n$-precluster tilting in terms of relative cotilting. In this section we give an alternate definition more closely linked to $n$-cluster tilting, namely that ${}^{\perp_{n-1}}\M = \M^{\perp_{n-1}}$. We begin by connecting this to relative cotilting by showing that ${}^{\perp_{n-1}}\M={}^{\perp_{F_\M}}\M$.

\begin{lemma}
  Let $n>1$, and let $\M$ be $n$-rigid, generating, and contravariantly finite. Then $\Ext^i(-, M) = \Ext^i_{F_\M}(-, M)$ for all $M \in \M$ and $i < n$. In particular, if $\M$ is $n$-precluster tilting we have ${}^{\perp_{n-1}}\M={}^{\perp_{F_\M}}\M$.
  \begin{proof}
    Let $X$ be an arbitrary object and let $\Omega_\M X \rightarrowtail M_X \twoheadrightarrow X$ be a right $\M$-approximation of $X$. Taking the long exact sequence in $\Ext(-,M)$ we see that $\Ext^1(X, M)$ is the cokernel of $$(M_X, M) \to (\Omega_\M X, M),$$ which means that $\Ext^1(X, M) = \Ext^1_{F_\M}(X, M)$. We also get 
    $$\Ext^{i-1}(\Omega_\M X, M) = \Ext^i(X, M)$$
    for all $1<i<n$, thus by induction $\Ext^i(X, M) = \Ext^i_{F_\M}(X, M)$ for all $i<n$.

    When $\M$ is $n$-precluster tilting, then by definition $\Ext^i_{F_\M}(X, M)=0$ for all $i\geq n$, hence we have ${}^{\perp_{n-1}}\M={}^{\perp_{F_\M}}\M$.
  \end{proof}
\end{lemma}

An important part of the definition of $n$-cluster tilting is that an $n$-cluster tilting subcategory $\M$ satisfies ${}^{\perp_{n-1}}\M = \M = \M^{\perp_{n-1}}$. This does not in general hold for $n$-precluster tilting subcategories, but we have a weaker condition.

\begin{proposition}
  If $\M$ is $n$-precluster tilting, then ${}^{\perp_{n-1}}\M = \M^{\perp_{n-1}}$.
  \begin{proof}
    First assume $X$ is in $\M^{\perp_{n-1}}$. Let $X \rightarrowtail M$ be an inflation into an object of $\M$. Then the cokernel is in $\M^{\perp_{n-2}}$. Iterating this $n-1$ times we get an $F_\M$-exact sequence
    \begin{center}
      \begin{tikzcd}
        X \ar[r, rightarrowtail] & M_1 \ar[r] & M_2 \ar[r] & \cdots \ar[r] & M_{n-2} \ar[r] & M_{n-1} \ar[r, twoheadrightarrow] & Y,
      \end{tikzcd}
    \end{center}
    thus $X = \Omega_{F_\M}^{n-1} Y$. Then $\Ext_{F_\M}^{>0}(X, M) = \Ext^{>n-1}_{F_\M}(Y, M)=0$, so $X$ is in ${}^{\perp_{F_{\M}}}\M = {}^{\perp_{n-1}}\M$.

    Next assume $X$ is in ${}^{\perp_{n-1}}\M = {}^{\perp_{F_{\M}}}\M$. Then we have a conflation $X \rightarrowtail M \twoheadrightarrow Y$ with $M \twoheadrightarrow Y$ a right $\M$-approximation, and $Y$ also in ${}^{\perp_{F_{\M}}}\M$. The long exact sequence in $\Ext(M,-)$ then give us that $\Ext^1(M, X)=0$ and $\Ext^{i}(M, X)=\Ext^{i-1}(M, Y)$ for $1 < i < n$. Thus by induction $X$ is in $\M^{\perp_{n-1}}$.
  \end{proof}
\end{proposition}

In fact this property is strong enough to give an alternative definition of precluster tilting, by replacing the cotilting condition. This extends the result of \cite[Prop~3.8]{IS18} to the exact setting.

\begin{proposition}\label{prop:n-precluster_tilt_iff_symmetric_orthogonality}
  Let $n>1$ and let $\M$ be a generating-cogenerating, functorially finite, $n$-rigid subcategory of an exact category $\E'$. Then ${}^{\perp_{n-1}}\M = \M^{\perp_{n-1}}$ if and only if $\M$ is $n$-precluster tilting.
  \begin{proof}
    We have already seen one direction, so assume ${}^{\perp_{n-1}}\M = \M^{\perp_{n-1}}$. Then we need to show that $\id_{F_{\M}} \M < n$, and that $\M$ is $F_\M$-cotilting. We first show that $\id_{F_{\M}} \M < n$.

    Let $X$ be an arbitrary object of $\E'$, and let $M$ be an object of $\M$. Since $\M$ is $n$-rigid we have that $\Ext^{i}(X, M) = \Ext^{i}_{F_\M}(X, M)$ for $i<n$. Taking an $\M$-approximation of $X$ we get an exact sequence 
    \begin{center}
      \begin{tikzcd}
        0 \ar[r] & \Omega_\M X \ar[r] & M_X \ar[r] & X \ar[r] & 0.
      \end{tikzcd}
    \end{center}
    Then from the long exact sequence in $\Ext(M, -)$ we get $\Ext^1(M, \Omega_\M X)=0$ and $\Ext^{i-1}(M, X)=\Ext^i(M, \Omega_\M X)$ for all $1<i<n$. Therefore, by induction, $\Omega^{n-1}_\M X$ is in $\M^{\perp_{n-1}} = {}^{\perp_{n-1}}\M$, and we have 
    $$0 = \Ext^1(\Omega^{n-1}_\M X, M) = \Ext^1_{F_{\M}}(\Omega^{n-1}_\M X, M) = \Ext^n_{F_{\M}}(X, M).$$
    Since $X$ is arbitrary this means that $\id_{F_\M} M < n$.

    Next, we show that $\M$ is $F_\M$-cotilting. Let $X$ be in ${}^{\perp_{F_{\M}}}\M \subseteq {}^{\perp_{n-1}}\M = \M^{\perp_{n-1}}$. If we take a left $\M$-approximation we get a conflation $X \rightarrowtail M \twoheadrightarrow Y$. Since $\Ext^1(M ,X)=0$ this conflation is $F_\M$-exact. Since $X \rightarrowtail M$ is a left $\M$-approximation we get that $Y$ is in ${}^{\perp_{F_{\M}}}\M$. Thus $\M$ is $F_\M$-cotilting, and consequently $n$-precluster tilting.
  \end{proof}
\end{proposition}

We conclude this section by proving an analogues result to \cite[Prop~3.12(b)]{IS18}.

\begin{proposition}
  Let $\M$ be $n$-precluster tilting. Consider ${}^{\perp_{F_{\M}}}\M$ as an exact category with the $F_\M$-exact structure. Then ${}^{\perp_{F_{\M}}}\M$ is Frobenius with projective objects $\M$.
  \begin{proof}
    Since $\M$ is generating and contravariantly finite ${}^{\perp_{F_{\M}}}\M$ has enough projectives, with $\M$ as the projective objects. Since $\M$ is $F_\M$-cotilting ${}^{\perp_{F_{\M}}}\M$ has enough injectives with $\M$ as the injective objects. Thus ${}^{\perp_{F_{\M}}}\M$ is Frobenius with projective objects $\M$.
  \end{proof}
\end{proposition}

\section{Comparisons to the classical case}\label{sec:comparison_to_classical_case}

In this section we compare our definition of precluster tilting and minimal Auslander--Gorenstein to those of Iyama--Solberg for Artin algebras\cite{IS18}. Throughout this section the symbols $\Gamma$ and $\Lambda$ refer to Artin algebras. We begin by recalling the definitions.

\begin{definition}[{}$n$-minimal Auslander--Gorenstein]\cite[Def~1.1]{IS18}
  We call an Artin algebra $\Gamma$ an $n$-minimal Auslander--Gorenstein algebra if it satisfies
  $$\id_\Gamma \Gamma \leq n+1 \leq \domdim \Gamma.$$
\end{definition}

\begin{definition}[{$n$}-precluster tilting]\cite[Def~3.2]{IS18}
  A subcategory $\mathcal C$ of $\mod \Lambda$ is called an $n$-precluster tilting subcategory if it satisfies the following conditions.
  \begin{enumerate}[label=(\roman*)]
    \item $\mathcal C$ is a generator-cogenerator for $\mod\Lambda$,
    \item $\tau_n(\mathcal C) \subseteq \mathcal C$ and $\tau_n^- (\mathcal C) \subseteq \mathcal C$,
    \item $\Ext^i(\mathcal C, \mathcal C) = 0$ for $0 < i < n$,
    \item $\mathcal C$ is a functorially finite subcategory of $\mod\Lambda$.
  \end{enumerate}
  If moreover $\mathcal C$ admits an additive generator $M$, we say that $M$ is an $n$-precluster tilting module.
\end{definition}

Notice that the only difference here is that the condition for $\mathcal C$ to be $F_{\mathcal C}$-cotilting is replaced by $\mathcal C$ being closed under the higher Auslander--Reiten translations, $\tau_n := \tau\Omega^{n-1}$ and $\tau_n^- := \tau^-\mho^{n-1}$. 

\subsection{The correspondence for Artin algebras}

One distinction between the higher Auslander--Solberg correspondence of Iyama--Solberg\cite{IS18} and ours is that Iyama--Solberg uses the contravariant Yoneda embedding to go from an $n$-precluster module to an $n$-minimal Auslander--Gorenstein algebra. This gives an equivalence between the $n$-precluster subcategory and the opposite category of projectives in the corresponding algebra. This is quite natural when working with modules since $\Hom(\Lambda, M) \cong M$, and in the presence of a duality the two approaches are equivalent, because $\Hom(M, D\Lambda) = DM$. However we use the covariant Yoneda embedding, building on \cite{henrard2020auslanders} and \cite{ebrahimi2021higher}. This causes some applications of duality when comparing to the Artin case (c.f. \cref{cor:artin_n-auslander-gorenstein}), but is not important since the definitions are symmetric over Artin algebras.

We begin by showing our defintion of Morita--Tachikawa categoris is an honest generalization of algebras with dominant dimension at least 2.

\begin{proposition}\label{prop:classical_morita-tachikawa}
  Let $\Gamma$ be an Artin algebra. Then $\Gamma$ has dominant dimension at least 2 if and only if $\mod\Gamma$ is Morita--Tachikawa.
  \begin{proof}
    If $\Gamma$ has dominant dimension at least 2, then $\mod\Gamma$ being Morita--Tachikawa follows from \cite[Lemma~4.14, Lemma~4.15, and Lemma~4.17(2)]{henrard2020auslanders}.
    
    Conversely, assume $\mod\Gamma$ is Morita--Tachikawa and let 
    \begin{center}
      \begin{tikzcd}
        0 \ar[r] & \Gamma \ar[r] & I_0 \ar[r] & I_1
      \end{tikzcd}
    \end{center}
    be a minimal injective copresentation. Since $({}^\perp \Gamma, \cogen\Gamma)$ is a torsion pair and the socle of $I_0$ is in $\cogen\Gamma$, we must have $I_0$ in $\cogen\Gamma$. Thus $I_0$ is projective. Now, using the torsion pair, let $X$ be the socle of the torsion part of $I_1$. Then 
    \begin{align*}
      \Hom(X, I_1) = \Ext^1(X, \Gamma) = 0
    \end{align*}
    and so $X=0$. Thus $I_1$ is in $\cogen\Gamma$ and $\Gamma$ has dominant dimension at least 2.
  \end{proof}
\end{proposition}

Next we show that the correspondence in \cref{sec:Auslander-Solber-correspondence} restricts to the classical Morita--Tachikawa correspondence (up to duality).

\begin{proposition}\label{prop:classical_correspondence}
  Let $\Gamma$ be an Artin algebra with dominant dimension at least 2. Let $I$ be an additive generator for the subcategory of projective-injective modules, and let $\Lambda$ be the endomorphism ring $\End_\Gamma(I)$. Then we have a commutative diagram of exact functors:
  \begin{center}
    \begin{tikzcd}
      \mod \Gamma \ar[swap]{dr}{D\Hom_\Gamma(-,I)} \ar[r] & \left.(\mod \Gamma) \middle/ {}^\perp \Gamma\right. \ar[d, "\cong"] \\
      & \mod \Lambda.
    \end{tikzcd}
  \end{center}
  \begin{proof}
    Since the injective envelope of $\Gamma$ is in $\add I$, we have ${}^\perp \Gamma \subseteq {}^\perp I$. Then, since $I$ is injective $D\Hom_\Gamma(-,I)$ is exact, and so we get an induced exact functor from $\left.(\mod \Gamma) \middle/ {}^\perp \Gamma\right.$ to $\mod \Lambda$ making the diagram commute. We just need to show that this functor is an equivalence.

    We saw in \cref{lem:equiv_on_syzygy} that the localization functor has a fully faithful right adjoint with essential image the subcategory of second syzygies. The functor $D\Hom_\Gamma(-,I)$ also has a fully faithful right adjoint, given by $\Hom_\Lambda(DI, -)$, whose essential image consists of those objects with a copresentation in $\add I$. If we can show that the two images coincides, then it follows that the induced functor is an equivalence.

    Any object with a copresentation in $\add I$ is a second syzygy, since $I$ is projective. Now, let $X$ be a second syzygy. Then $X$ fits into an exact sequence 
    \begin{center}
      \begin{tikzcd}
        0 \ar[r] & X \ar[r] & P \ar[r] & P'
      \end{tikzcd}
    \end{center}
    with $P$ and $P'$ projective. Since all projectives embeds into an object in $\add I$ we may assume $P'$ is in $\add I$. Let $I(P)$ be the injective envelope of $P$, and take a pushout 
    \begin{center}
      \begin{tikzcd}
        0 \ar[r] & X \ar[r] \ar[d, equal] & P \ar[r] \ar[d] & P' \ar[d]\\
        0 \ar[r] & X \ar[r] & I(P) \ar[r] & PO \arrow[ul, phantom, "\scalebox{1}{$\ulcorner$}" , very near start, color=black]
      \end{tikzcd}
    \end{center}
    Since $PO$ is the cokernel of $P \to I(P) \oplus P'$ and $\Gamma$ has dominant dimension at least 2, we have that $PO$ embeds into an object of $\add I$. Thus $X$ has a copresentation in $\add I$, and the induced functor is an equivalence.
  \end{proof}
\end{proposition}

The last puzzle piece we need is just a slight generalization of what we already proved in \cref{prop:n-precluster_tilt_iff_symmetric_orthogonality}.

\begin{proposition}
  If $\M$ is a subcategory of $\mod\Lambda$ for an Artin algebra $\Lambda$, then $\M$ is an $n$-precluster tilting subcategory if and only if it satisfies \cite[Def~3.2]{IS18}.
  \begin{proof}
    When $n>1$ this follows from \cref{prop:n-precluster_tilt_iff_symmetric_orthogonality} and \cite[Prop~3.8(b)]{IS18}, so we need only consider the case $n=1$.

    Assume $\M$ is functorially finite and generating-cogenrating. Then $\M$ is $1$-precluster tilting if and only if $\M$ is presicely the $F_\M$-injectives and $(\mod\Lambda, F_\M)$ has enough injectives. The $F_{\M}$-injectives are given by $\add(\tau\M \cup \add D\Lambda)$ \cite[Prop~1.9]{AS93}, so if $\M$ is $1$-precluster tilting, then $\tau \M \subseteq \M$ and $\tau^- \M \subseteq \M$. For the converse, if $\tau \M \subseteq \M$ and $\tau^- \M \subseteq \M$, then $\M$ is equal to $\add(\tau\M \cup \add D\Lambda)$, which is the class of $F_\M$-injectives. By \cite[Prop~1.8, Prop~1.12(b)]{AS93}, since $\M$ is covariantly finite $(\mod\Lambda, F_\M)$ has enough injectives.
  \end{proof}
\end{proposition}

Because of the correspondence we have set up, this gives us the analogues statement for $n$-minimal Auslander--Gorenstein algebras for free. We state it here as a corollary.

\begin{corollary}\label{cor:artin_n-auslander-gorenstein}
  If $\Gamma$ is an Artin algebra, then $\Gamma$ is $n$-minimal Auslander--Gorenstein if and only if $\mod\Gamma$ is an $n$-minimal Auslander--Gorenstein category.
  \begin{proof}
    We saw in \cref{prop:classical_morita-tachikawa} that $\Gamma$ has dominant dimension 2 if and only if $\mod\Gamma$ is Morita--Tachikawa, so we may assume this. Then \cref{prop:classical_correspondence} and \cref{thm:Auslander-Solberg_correspondence} gives us that $\mod\Gamma$ is $n$-minimal Auslander--Gorenstein if and only if $\add DI$ is $n$-precluster tilting in $\mod\End(I)$, where $I$ is an additive generator for the subcatgeory of projective-injective $\Gamma$-modules. At the same time, \cite[Thm~4.5]{IS18} tells us that $DI$ is $n$-precluster tilting if and only if $\Gamma$ is $n$-minimal Auslander--Gorenstein. Thus the two conditions are equivalent.
  \end{proof} 
\end{corollary}

\section{Injectives and projectives}\label{sec:projectives_and_injectives}

In this section we investigate the projective and the injective objects of $\E/\eff$ for a \pac{} $\E$. The results and proofs are almost identical to those in \cite[Section~4.3--4.4]{henrard2020auslanders}.

We begin by considering injective objects.

\begin{proposition}\cite[Lemma~4.11]{henrard2020auslanders}\label{prop:projective_injectives_in_auslander_category}
  Let $\E$ be an exact category and let $\M \subseteq \E$ be a generating-cogenerating contravariantly finite subcategory. Then $(-,I)$ is injective in $\mod_{\adm}(\M)$ if and only if $I$ is injective in $\E$. 
  \begin{proof}
    Assume $I$ is injective, and let $F$ be in $\mod_{\adm}(\M)$. Then $F$ fits into an exact sequence 
    \begin{center}
      \begin{tikzcd}
        0 \ar[r] & (-,K) \ar[r] & (-, M') \ar[r] & (-,M) \ar[r] & F \ar[r] & 0.
      \end{tikzcd}
    \end{center}
    Applying $(-, (-,I))$ and using \cref{lemma:yoneda} we get
    \begin{center}
      \begin{tikzcd}
        0 \ar[r] & (F, (-,I)) \ar[r] & (M, I) \ar[r] & (M', I) \ar[r] & (K, I) \ar[r] & 0.
      \end{tikzcd}
    \end{center}
    Since $I$ is injective this sequence is exact and thus $\Ext^1(F, (-,I))=0$. Then it follows that $(-, I)$ is injective.

    Conversely, assume $I$ is not injective. Then there exists an inflation $I \rightarrowtail X$ that is not split. Since $\M$ is cogenerating this means there exists an inflation $I \rightarrowtail M$ with $M \in \M$ that is not split. This gives rise to an inflation $(-, I) \rightarrowtail (-, M)$ in $\mod_{\adm}(\M)$. By \cref{lemma:yoneda} this inflation is not split and so $(-,I)$ is not injective.
  \end{proof}
\end{proposition}

Next, we investigate the projective objects.

\begin{proposition}\cite[Lemma~4.25]{henrard2020auslanders}\label{prop:projectives_in_category_with_n-precluster}
  Let $\E$ be a \pac{} with projective objects $\P$. Then $P\in \P$ is projective in $\E/\eff$ if and only if $P$ is in ${}^\perp (\eff)$.
  \begin{proof}
    Assume that $P \in \P$ is not projective in $\E/\eff$. Then there is some nonsplit deflation $X \twoheadrightarrow P$ in $\E/\eff$. By \cref{lem:equiv_on_syzygy} this corresponds to an exact sequence $X \to P \to C \to 0$ in $\E$ with $C \in \eff$. Since $P$ is projective in $\E$ and the map is not split, $C$ must be nonzero. Therefore $P$ is not in ${}^\perp (\eff)$.

    Conversely, if $P$ is not in ${}^\perp (\eff)$ then there is some non-zero map $P \to C$ with $C\in\eff$. Let 
    \begin{center}
      \begin{tikzcd}
        P_1 \ar[r]& P_0 \ar[r] & C \ar[r] & 0
      \end{tikzcd}
    \end{center}
    be a projective presentation of $C$ in $\E$. Then $P_1 \to P_0$ is a deflation in $\E/\eff$. Since $\Hom_\E(P, C) \neq 0$, the map 
    \begin{center}
      \begin{tikzcd}
        \Hom_\E(P, P_1) \ar[r] & \Hom_\E(P, P_0)
      \end{tikzcd}
    \end{center}
    is not surjevctive. By \cref{lem:equiv_on_syzygy} we have $\Hom_\E(P, P_i) = \Hom_{\E/\eff}(P, P_i)$, so since $P_1 \to P_0$ is a deflation and the above map is not surjective, it follows that $P$ is not projective in $\E/\eff$.
  \end{proof}
\end{proposition}

Now that we have an understanding of what the projective and the injective objects are we bring our attention to when we have enough projectives/injectives.

\begin{proposition}\cite[Lemma~4.13]{henrard2020auslanders}
  Let $\M$ be a generating-cogenerating contravariantly finite subcategory of an exact category $\E'$. Then $\E'$ has enough injectives if and only if $\domdim \mod_{\adm}(\M) \geq 1$.
  \begin{proof}
    Assume $\E'$ has enough injectives. Then for every object $M \in \M$ we have an inflation $M \rightarrowtail I$ for an injective object $I$, which induces an inflation $(-,M) \rightarrowtail (-,I)$. By \cref{prop:projective_injectives_in_auslander_category} $(-,I)$ is projective-injective. Since all projectives in $\mod_{\adm}(\M)$ are of the form $(-,M)$ for some $M\in\M$ we get that $\domdim \mod_{\adm}(\M) \geq 1$.

    Conversely, since $\M$ is cogenerating any object of $\E$ embeds into an object of $\M$. Since $\domdim \mod_{\adm}(\M) \geq 1$ each object of $\M$ embeds into an injective, and therefore so does any object.
  \end{proof}
\end{proposition}

\begin{proposition}\cite[Prop~4.26]{henrard2020auslanders}
  Let $\E$ be a \pac{} with projectie objects $\P$. Let $\mathcal G$ denote the subcategory $\operatorname{gen}\left(\P\cap {}^\perp (\eff))\right) \subseteq \E$. Then $\E/\eff$ has enough projectives if and only if $\left(\mathcal G, \eff \right)$ forms a torsion pair in $\E$.
  \begin{proof}
    First assume $\E/\eff$ has enough projectives. It's clear that $\Hom_\E(\mathcal G, \eff) =0$, so we just need to prove that every object $X\in\E$ fits into an exact sequence 
    \begin{center}
      \begin{tikzcd}
        0 \ar[r] & G \ar[r] & X \ar[r] & C \ar[r] & 0
      \end{tikzcd}
    \end{center}
    with $G$ in $\mathcal G$ and $C$ in $\eff$. Let $X \xrightarrow{\sim} F$ be a weak equivalence to a second syzygy as in \cref{lem:equiv_on_syzygy}. Since $\E/\eff$ has enough projectives and the projectives are given by $\P\cap {}^\perp (\eff)$ we can find an object $P \in \P\cap {}^\perp (\eff)$ and a map $P \to F$ that becomes a deflation in $\E/\eff$. Since the cokernel of $X \xrightarrow{\sim} F$ is in $\eff$ and $P$ is in ${}^\perp (\eff)$ the map $P\to F$ factors thorugh the image of the weak equivalence. Since $P$ is projective $P\to F$ factors through $X$. Lastly, since the localization reflects admissible morphisms \cite[Thm~2.16]{henrard2020auslanders}, the map $P\to X$ is admissible with cokernel in $\eff$. Letting $G$ be the image of $P\to X$ we get our desired sequence.

    Conversely, if $\left(\mathcal G, \eff \right)$ forms a torsion pair then $\mathcal G$ becomes dense in $\E/\eff$, and so $\E/\eff$ has enough projectives.
  \end{proof}
\end{proposition}

\subsection[Flavors of \texorpdfstring{$n$}{n}-minimal Auslander--Gorenstein categories]{Flavors of \texorpdfstring{\boldmath $n$}{n}-minimal Auslander--Gorenstein categories}

In \cite[Section~6]{IS18} Iyama and Solberg classify $n$-precluster tilting modules into 4 distinct classes, and describe the corresponding classes of $n$-minimal Auslander--Gorenstein algebras. Specifically they say that for an Artin algebra $\Lambda$, an $n$-precluster tilting module $M$ satisfies $\add M = \add \{\P_n \wedge \mathcal I_n, N\}$ where $\P_n \wedge \mathcal I_n$ is the unique minimal $n$-precluster tilting module, and $N$ some module satisfying $\tau_n(N) \cong N$. Then one of four distinct cases can occur:
\begin{enumerate}[label=(\Alph*)]
  \item $\Lambda$ is selfinjective and $N=0$.
  \item $\Lambda$ is selfinjective and $N\neq 0$.
  \item $\Lambda$ is not selfinjective and $N=0$.
  \item $\Lambda$ is not selfinjective and $N\neq0$.
\end{enumerate}

Flavor A has an obvious generalization to the exact case. Namely, that for a Frobenius catgeory, the subcatgeory of projectives is the unique minimal $n$-precluster tilting subcategory for any $n$. In this case $\mod_{\adm}(\M)$ is just equivalent to the original category.

For flavor B, there are a few possibilites for the `correct' generalization. The most obvious is perhaps Frobenius categories that are not of flavor A, but we propose a slightly more general defintion. Namely, that $(\E', \M)$ is of flavor $B$ if it is not of flavor A and the injective and projective modules of $\E'$ coincide. From \cref{prop:projective_injectives_in_auslander_category} and \cref{prop:projectives_in_category_with_n-precluster} we see that this is equivalent to $\mod_{\adm}(\M)$ having ${}^\perp (\eff) \cap \P$ as projective-injective objects without being Frobenius.

For flavor C and D it is not clear to the author what the correct generalization should be. We do not suspect that a unique minimal $n$-precluster tilting subcategory should exist in general, and it is not clear what it would say about $\mod_{\adm}(\M)$ if it did. We frame this as an open problem.

\begin{problem}
  If $\E'$ is a category that admits an $n$-precluster tilting subcategory. Does $\E'$ admit a minimal $n$-prelcuster tilting subcatgeory, and under what conditions is this subcategory unique?
\end{problem}

\section{Examples}\label{sec:examples}

In this section we give some explicit examples of $n$-precluster tilting subcategories and $n$-minimal Auslander--Gorenstein categories. We first restate an example from the previous section.

\begin{example}
  Any Frobenius category is $n$-minimal Auslander--Gorenstein, and the subcategory of projective objects is $n$-precluster tilting for all $n$. In this case the bijection of \cref{thm:bijection_partial_auslander} is trivial, since $\eff = 0$.
\end{example}

Another example we have seen before is that $\mod \Gamma$ is an $n$-minimal Auslander--Gorenstein exact category whenever $\Gamma$ is an $n$-minimal Auslander--Gorenstein algebra, and that $\add M$ is an $n$-precluster tilting subcategory whenever $M$ is an $n$-precluster tilting module. We give a specific example of this below.

\begin{example}
  Let $\Lambda:=\Pi A_3$ be the preprojective algebra of $1\to2\to3$, and let $M$ be the (left) $\Lambda$-module $M:=\Lambda \oplus P_2/S_2$. Then $M$ is 2-precluster tilting.

  Consequently, the algebra $\Gamma:=\End(M)$ is 2-Auslander--Gorenstein given by the quiver
  \begin{center}
    \begin{tikzcd}
      1 \ar[rd, swap, "\alpha"] & 2 \ar[l, swap, "\alpha^*"] \ar[r, "\beta"] & 3 \ar[dl, "\beta^*"]\\
      & 4 \ar[u, swap, "\gamma"]
    \end{tikzcd}
  \end{center}
  with relations $\beta\gamma\beta^*$, $\alpha^*\gamma\alpha$, $\alpha\alpha^*-\beta^*\beta$ and $\alpha\alpha^*\gamma$.
  
  Notice that $\Gamma$ has infinite global dimension (since $\Omega^2 S_1 = S_1$), so $\Gamma$ is not 2-Auslander. Therefore $M$ is not 2-cluster tilting, which we can also see by noting that $$\Ext^1_\Lambda(M, S_2) = \Ext^1_\Lambda(S_2, M) = 0.$$
\end{example}

Continuing with this example, let $\Gamma$ be the Gorenstein-algebra we just defined. The Gorenstein projectives of $\Gamma$ are the projective modules as well as $S_1$, $S_3$, $JP_1$, and $JP_3$. Since this forms a Frobenius category the category of projectives is $n$-precluster tilting for all $n$, but there are other $n$-precluster tilting subcategories.
\begin{example}
  The following are $1$-precluster tilting subcategories of the category of Gorenstein projectives over $\Gamma$:
  \begin{center}
    \renewcommand{\labelitemi}{--}
    \begin{itemize}
      \item $\add(\Gamma \oplus S_1 \oplus JP_1 )$
    \item $\add(\Gamma \oplus S_3 \oplus JP_3 )$
    \end{itemize}
    \end{center}
  The following are $2$-precluster tilting:
  \begin{multicols}{2}
    \begin{center}
    \renewcommand{\labelitemi}{--}
    \begin{itemize}
      \item $\add(\Gamma \oplus S_1)$    
      \item $\add(\Gamma \oplus JP_1)$
      \item $\add(\Gamma \oplus S_1 \oplus S_3)$ 
      \item $\add(\Gamma \oplus S_3)$
      \item $\add(\Gamma \oplus JP_3)$
      \item $\add(\Gamma \oplus JP_1\oplus JP_3)$
    \end{itemize}
    \end{center}
  \end{multicols}
  The bottom two are also $2$-cluster tilting.
\end{example}

In \cite{Iyama11}, Iyama gives several exmaples of cluster tilting modules in $T^{\perp_\infty}$ for a tilting module $T$, which they call relative cluster tilting. In a private conversation Sondre Kvamme suggested an idea of how to use this to construct `relative' precluster tilting. The idea is that if we start with a relative cluster tilting module and take the tensor with a selfinjective algebra, we should get a precluster tilting module in the $\Ext$-complement of a tilting module. We give a proof of this below.

\begin{proposition}
  Let $\Lambda$ be a finite dimensional algebra over a field $k$, and let $T$ be a tilting module over $\Lambda$. Let $M$ be $n$-precluster tilting in $T^{\perp_\infty}$ for some $n>1$, and let $\Gamma$ be a selfinjective algebra over $k$. Then $M\otimes\Gamma$ is $n$-precluster tilting in $(T\otimes\Gamma)^{\perp_\infty}$.
  \begin{proof}
    Let $X$ be any $\Lambda$-module and let $Y$ be any $\Lambda\otimes\Gamma$-module. Then we have 
    $$X\otimes\Gamma = X \mathop\otimes\limits_\Lambda (\Lambda \otimes \Gamma).$$ 
    If we use this together with the Hom-Tensor adjunction, it gives us 
    $$ \Hom_{\Lambda \otimes \Gamma}(X \otimes \Gamma, Y) = \Hom_{\Lambda}(X, Y). $$
    Using duality, and the fact that $\Gamma \cong D\Gamma$ we also get
    \begin{align*}
      \Hom_{\Lambda \otimes \Gamma}(Y, X \otimes \Gamma) &= \Hom_{\Lambda \otimes \Gamma}(DX \otimes D\Gamma, DY) \\
      &= \Hom_{\Lambda}(DX, DY) \\
      &= \Hom_\Lambda(Y, X).
    \end{align*}
    With these two facts one can show that $T\otimes\Gamma$ is a tilting module, and that $M\otimes\Gamma$ is $n$-rigid.
  
    Now, since whenever $P_M^\bullet$ is a projective resolution of $M$, we have that $P_M^\bullet \otimes\Gamma$ is a projective resolution of $M\otimes\Gamma$, we get
    $$ \Ext_{\Lambda\otimes\Gamma}^i(M\otimes\Gamma, Y) = \Ext_\Lambda^i(M, Y)
      \quad\text{and}\quad
      \Ext_{\Lambda\otimes\Gamma}^i(Y, M\otimes\Gamma) = \Ext_\Lambda^i(Y, M).$$
    Thus ${}^{\perp_{n-1}}M = M^{\perp_{n-1}}$ also in $(T\otimes\Gamma)^{\perp_\infty}$, and we conclude that $M\otimes\Gamma$ is $n$-precluster tilting.
  \end{proof}

\end{proposition}

\clearpage

\bibliography{exactbib}
\bibliographystyle{alpha}

\end{document}